\documentclass[11pt]{article}
\usepackage{amsfonts}
\usepackage{bookmark}

\usepackage{amsmath,amsfonts,amssymb,amscd,amsthm,amsbsy,bm,latexsym}
\usepackage [latin1]{inputenc}
\usepackage[pagewise]{lineno}

\newtheorem{lemma}{Lemma}[section]
\newtheorem{theorem}{Theorem}[section]

\newtheorem{corollary}{Corollary}[section]

\numberwithin{equation}{section}

\begin{document}

\title{Blow-up phenomena for a class of extensible beam equations   \footnote{ Email addresses:  gongweiliu@haut.edu.cn}}

\author
{Gongwei Liu, \,\,Mengyun Yin,\,\,Suxia Xia\\
School of Science, Henan University of Technology, Zhengzhou 450001, China}

\date{}
\vskip 0.2cm
\maketitle
\noindent{\bf Abstract}: In this paper, we  investigate the  initial boundary value problem of the following  nonlinear extensible beam equation with nonlinear damping term 
$$u_{t t}+\Delta^2 u-M\left(\|\nabla u\|^2\right) \Delta u-\Delta u_t+\left|u_t\right|^{r-1} u_t=|u|^{p-1} u$$
which was considered by Yang et al. (Advanced Nonlinear Studies 2022; 22:436-468). We consider the problem with the nonlinear damping and  establish the finite time blow-up of the solution for the initial data at arbitrary high energy level, including   the estimate lower and upper bounds of the blowup time.  The result  provides some affirmative  answer to the open problems given in (Advanced Nonlinear Studies 2022; 22:436-468).



\section{Introduction}
\setcounter{equation}{0}
In this paper, we are concerned with  the initial boundary value problem of the following beam equation with linear strong damping, nonlinear weak damping and nonlinear source
\begin{equation}\label{1.1}
\begin{cases}
u_{t t}+\Delta^2 u-M\left(\|\nabla u\|^2\right) \Delta u-\Delta u_t+\left|u_t\right|^{r-1} u_t=|u|^{p-1} u,&(x,t)\in \Omega_T, \\
u(x, 0)=u_0(x),\,\,\, u_t(x, 0)=u_1(x),\,\,\qquad \qquad &x\in \Omega,\\
u=\frac{\partial u}{\partial\nu}=0, \qquad \qquad \qquad\quad&(x, t)\in\partial\Omega,
\end{cases}
\end{equation}
where $\Omega\subset\mathbb{R}^n$ ($n\geq1$) is a bounded domain with  smooth boundary $\partial \Omega$, $\Omega_T:=\Omega\times(0,T)$, the nonlinear Kirchhoff term $M(s)$ is defined by
\begin{equation}\label{1.2}
M(s)=1+\beta s^\gamma, \gamma \geq 0, \beta \geq 0, s \geq 0,
\end{equation}
and the exponent $p$ of the source term satisfies
 \begin{equation}\label{1.3}
1\leq r<p,\quad 1<2\gamma+1<p<\infty, \quad n\leq4 ; \quad 1<2 \gamma+1<p \leq \frac{n+2}{n-4}, \quad n \geq 5.
\end{equation}
The motivation for the problem \eqref{1.1} comes from the following original equation $$
\frac{\partial^2 u}{\partial t^2}+\alpha \frac{\partial^4 u}{\partial x^4}-\left(\beta+k \int_0^L u_x^2 \mathrm{~d} x\right) \frac{\partial^2 u}{\partial x^2}=0,
$$
 proposed by Woinowsky-Krieger \cite{WK1950} in 1950. Hence, the model \eqref{1.1} can be seen as the nonlinear extensible beam equation in a type of modified Woinowsky-Krieger.
For more physical background of \eqref{1.1} and the related models, we can refer the readers to \cite{BallJMAA1973,BallJDE1973,Berger,CT2005,YangANS2022} and the references therein. It is well-known that there are lots of interesting results ( such as local/global existence of the solutions, blow-up and the lifespan,  energy decay, longtime dynamics and so on under suitable conditions) on the extensible beam models with different damping term (such as linear or nonlinear damping, strong damping). Here, we only list partial results \cite{CaDo2004,ChenZhouNA2009,GuoTWJ2019,HanEECT2022,LiuERA2020,Mess2002,WuTWJ2009,YangJDE2013}. For more details, we  refer
the readers to the Introduction section in \cite{YangANS2022}.

 The present paper mainly deals with the effect of the initial data on the finite time blow up of the solution to the problem \eqref{1.1} with nonlinear damping term under the arbitrarily high initial energy. The main goal of the present paper is to give some affirmative  answer to the open problem given in \cite{YangANS2022}, so we just make a quick start from \cite{YangANS2022}. In the frame work of potential well, the authors of \cite{YangANS2022} established  global existence, nonexistence, and asymptotic behavior of solution for both subcritical and critical initial energy level. Moreover, they also obtained the global
nonexistence for the problem with linear weak damping ($r=1$) and strong damping at the
high initial energy level. However, for  the high energy case, i.e., $E(0)>0$, the global nonexistence of the problem with  nonlinear damping term $|u_t|^{r-1}u_t, r>1$ is still an open problem \cite{YangANS2022}. Moreover, there is no information on upper or lower bound of the blow-up time.
The main purpose of this paper is to give a positive answer to this open problem and  generalize some results obtained  in \cite{YangANS2022}.

As far as we know, for the wave-type equations with nonlinear damping, there are few results on the finite time time blow up and the  bounds of the blowup time  under the arbitrarily high initial energy.
 Recently, Sun et al \cite{SunAA2019} investigated  a  nonlinear viscoelastic wave equation with nonlinear damping and obtained blow-up of solutions when the initial
data at  high energy level. Similar results were also obtained for Petrovsky
type equation by Liu et al \cite{LiuBOVP2019}. It is also worthy pointing out that the recent  literatures \cite{LiZhouMJM2023} and \cite{LiaoLiMJM2023} where the blow-up of solutions to two types of  fourth-order
equation with variable-exponent nonlinear damping  and source   under high initial energy level were established. Moreover, Liao and Tan \cite{LiaoSCM2023} also got a lower bound of the blow-up time by making full use of the strong damping term $-\triangle u_t$.

Motivated by the above-mentioned literatures, in this research, we will establish some blow-up  results for problem \eqref{1.1} with arbitrary initial energy level (including $E(0)<0$ and $E(0)>d$). We also obtain an upper bound of the blow-up time under suitable condition.  Moreover, two lower bounds of the blow-up time are also obtained by different strategies.

 The rest of this paper can be organized  as follows. In Sect. 2, we give some notations and partial known results given in \cite{YangANS2022}, which will be used in this paper. In Sect. 3, we  establish the main results and proofs.

\section{ Preliminaries}
For  convenience, in this section,  we  collect some notations, functionals, assumptions and partial  results obtained in \cite{YangANS2022}.
Let   $\|\cdot\|_p$ indicate the norm in $L^p(\Omega)$ and  $(\cdot, \cdot)$ express the inner product in $L^2(\Omega)$. We use  the following   notations
\begin{equation*}
  H=\{u\in H^2(\Omega)\cap H^1_0(\Omega)|u=\frac{\partial u}{\partial \nu}=0\, \text{on}\,\,\partial \Omega\},\quad \|u\|_H^2=\|\nabla u \|_2^2+\|\triangle u\|_2^2.
\end{equation*}
When there is no possibility of confusion, we
denote by $C$ or $C_i$ a generic constant that may vary from line to line even in the same formula.
Denoting by $\lambda_1$ the first eigenvalue of the bi-harmonic operator with boundary condition $u|_{\partial\Omega}=\frac{\partial u}{\partial \nu}|_{\partial\Omega}=0$, then we have the following  inequalities \cite{EECT2017}
\begin{equation}\label{2.1}
 \lambda_{1}\|u\|^{2}_{2}\leqslant\| \triangle u\|^{2}_{2}, ~\lambda^{\frac{1}{2}}_{1}\|\nabla u\|^{2}_{2}\leqslant\|\triangle u\|^{2}_{2},~\forall u\in H.
\end{equation}

 We define the following potential functional and Nehari functional  as usual
\begin{equation}\label{2.2}
  J(u)=\frac{1}{2}\|\nabla u\|_2^2+\frac{1}{2}\|\triangle u\|_2^2+\frac{\beta}{2(\gamma+1)}\|\nabla u\|_2^{2\gamma +2}-\frac{1}{p+1}\|u\|_{p+1}^{p+1},
\end{equation}
\begin{equation}\label{2.3}
  I(u)=\|\nabla u\|_2^2+\|\triangle u\|_2^2+\beta\|\nabla u\|_2^{2\gamma +2}-\|u\|_{p+1}^{p+1}.
\end{equation}
We introduce the total energy functional for the problem \eqref{1.1}
\begin{equation}\label{2.4}
\begin{split}
 E(t)=&\frac{1}{2}\|u_t\|_2^2+\frac{1}{2}\|\nabla u\|_2^2+\frac{1}{2}\|\triangle u\|_2^2+\frac{\beta}{2(\gamma+1)}\|\nabla u\|_2^{2\gamma +2}-\frac{1}{p+1}\|u\|_{p+1}^{p+1}\\
 =& \frac{1}{2}\|u_t\|_2^2+\frac{p-1}{2(p+1)}(\|\nabla u\|_2^2+\|\triangle u\|_2^2)+\big(\frac{1}{2(\gamma+1)}-\frac{1}{p+1}\big)\beta\|\nabla u\|_2^{2\gamma +2}+\frac{1}{p+1}I(u).
\end{split}\end{equation}
By a  direct computation, we derive that
\begin{equation}\label{2.5}
E'(t)=-\|u_t\|_{r+1}^{r+1}-\|\nabla u_t\|_2^2\leq0.
\end{equation}
By $I(u)$, we define the stable set $W=\{u\in H|I(u)>0\}\cup \{0\}$ and the unstable set $V=\{u\in | I(u)<0\}$, respectively.  The depth of the potential well is defined by
\begin{equation*}
  d:=\inf_{u\in \mathcal{N}}J(u)=\inf_{u\in H\backslash \{\{0\}}\sup_{\lambda\geq0}J(\lambda u)
\end{equation*}
where $\mathcal{N}$ is the Nehari manifold  $\mathcal{N}=\{u\in H \backslash \{0\}|I(u)=0\}$. It follows from Lemma 2.3 in \cite{YangANS2022} that
$d=\frac{p-1}{2(p+1)}\big(\frac{1}{C^{p+1}}\big)^{\frac{2}{p-1}}$, where $C$ is the best embedding constant from $H$ into $L^{p+1}(\Omega)$, that is
\begin{equation}\label{2.6}
  \|u\|_{p+1}\leq C \|u\|_{H}.
\end{equation}

Now we state the following lemma.
\begin{lemma} Provided that $J(u)\leq d$, then $I(u)<0$ if and only if $\|u\|_H>\lambda_*:=C^{-\frac{p-1}{p+1}}$.
\end{lemma}
\begin{proof} If $\|u\|_H>\lambda_*:=C^{-\frac{p-1}{p+1}}$, \eqref{2.2}, \eqref{2.3} and the definition of $d$ yield that
\begin{equation*}
\begin{split}
J(u)&=\frac{p-1}{2(p+1)}(\|\nabla u\|_2^2+\|\triangle u\|_2^2)+\big(\frac{1}{2(\gamma+1)}-\frac{1}{p+1}\big)\beta\|\nabla u\|_2^{2\gamma +2}+\frac{1}{p+1}I(u)\\[2mm]
&\leq d=\frac{p-1}{2(p+1)}\big(\frac{1}{C^{p+1}}\big)^{\frac{2}{p-1}}=\frac{p-1}{2(p+1)}\lambda_*^2 .
\end{split}
\end{equation*}
Hence we can obtain $I(u)<0$ from the fact $p>2\gamma +1$ and $\beta\geq 0$.

On the other hand, assuming that $I(u)<0$, we have
\begin{equation*}
  \|u\|_H^2+\beta \|\nabla u\|_2^{2(\gamma+1)}<\|u\|_{p+1}^{p+1}\leq C^{p+1}\|u\|_{H}^{p+1},
\end{equation*}
which implies that $\|u\|_H>C^{-\frac{p-1}{p+1}}=\lambda_*.$ The proof is complete.
\end{proof}

To process our main results, we first state some results established  in \cite{YangANS2022} which will be used in the next section.
\begin{theorem}(\cite{YangANS2022} Local existence). Let $u_0(x)\in H$ and $u_1(x)\in H_0^1(\Omega)$. Then, the problem \eqref{1.1} admits a unique local solution $u:=u(x,t)\in C\big([0, T], H\big)$ satisfying  $u_t\in C\big([0,T]; H_0^1(\Omega)\big)\cap L^{r+1}\big([0, T], L^{r+1}(\Omega)\big)$ for some $T>0.$
\end{theorem}
Now combining Theorem 4.3 and Theorem 5.3 in \cite{YangANS2022}, the global nonexistence when $E(0)\leq d$ can be stated as follow.
\begin{theorem}(\cite{YangANS2022} Global nonexistence when $E(0)\leq d$). Let $u_0(x)\in H$ and $u_1(x)\in H_0^1(\Omega)$ be given functions. Assume that $E(0)\leq d$ and $u_0\in V$. Then, the solution $u$  to the problem \eqref{1.1} blows up in finite time.
\end{theorem}
By Lemma 2.1, we can derive that Theorem 2.2  is equivalent to the following form.
\begin{corollary}Provided that  $u_0(x)\in H$, $u_1(x)\in H_0^1(\Omega)$, $E(0)\leq d$ and $\|u_0\|_H> C^{-\frac{p-1}{p+1}}$, then  the solution $u$  to the problem \eqref{1.1} blows up in finite time.
\end{corollary}

 The global nonexistence of solution to the problem \eqref{1.1} with strong and linear damping ($r=1$) at the high  initial energy level $E(0)>0$ are obtained  by  invariant sign of $I(u(t))-\|u_t\|_2^2$.
\begin{theorem}(\cite{YangANS2022} Global nonexistence when $E(0)>0$ and $r=1$). Let $u_0(x)\in H$ and $u_1(x)\in H_0^1(\Omega)$ hold. Assume that $E(0)>0$, $I(u_0)-\|u_0\|_2^2<0$ and $\|\nabla u_0\|_2^2+\|u_0\|_2^2+2(u_0, u_1)>\frac{4(p+1)}{(p-1)\tilde{c}}$, where $\tilde{c}=\min\{1, C\}$ and $C$ is the best embedding constant from $H$ into $H_0^1(\Omega)$. Then, the solution $u$ to the problem \eqref{1.1}  with $r=1$ blows up in finite time.
\end{theorem}

\section{ Main Results and Proofs}
Now, we  state our main blow-up results as follows. Let $u(t)$ be the solution obtained in Theorem 2.1, whose maximal existence time is $T_m$.
\begin{theorem} Let \eqref{1.2} and \eqref{1.3} hold. Assume the initial data  $u_0(x)\in H$ and $u_1(x)\in H_0^1(\Omega)$ satisfies one of the following conditions
\begin{description}
  \item[(i)] $E((0)<0$;
  \item[(ii)] $0\leq E(0)<\frac{1}{B}\int_{\Omega}u_0u_1dx$.
\end{description}
Here $B$ is a  positive constant given in \eqref{3.1}. Then the solution $u$ to the problem \eqref{1.1} blows up in finite time.
\end{theorem}
\begin{proof}(i)  The  blow-up result for the case $E(0)<0$ is a direct conclusion of Theorem 2.2. Indeed, $E(0)<0$ and\eqref{2.4} yields that $I(u_0)<0$, i.e. $u_0\in V$.

(ii)  First, we assume that the energy $E(t)\geq0$ for all $t\in[0, T_m)$. Otherwise, there exist a $t_0\in[0, T_m)$ such that $E(t_0)<0$, taking $t_0$ as the initial time, by Case(i), we have that $u(t)$ blows up in finite time, which is a
contradiction.

We split the proof of (ii) into the following two steps.

\textbf{Step 1.} We first present the following claim which  is motivated from \cite{SunAA2019}.

\textbf{\emph{Claim}}: Suppose that $u_0(x)\in H$ and $u_1(x)\in H_0^1(\Omega)$ holds, and $u(t)$ is a weak solution to the  problem \eqref{1.1}. We
claim that there exist positive constants $A$ and $B$ such that
\begin{equation}\label{3.1}
  \frac{d}{dt}\big(\int_{\Omega}uu_tdx-BE(t)\big)\geq A \big(\int_{\Omega}uu_tdx-BE(t)\big) \text{ for all} \quad  t\in [0, T_m).
\end{equation}
\emph{Proof of the claim.} It follows from the first equation of the problem \eqref{1.1} that
\begin{equation*}
  \begin{split}\frac{d}{dt}\int_{\Omega}uu_tdx=&\|u_t\|_2^2+\int_{\Omega}uu_{tt}dx\\[2mm]
  =&\|u_t\|_2^2-\|\triangle u\|_2^2-\|\nabla u\|_2^2-\beta\|\nabla u\|_2^{2(\gamma+1)}-\int_{\Omega}\nabla u\cdot \nabla u_tdx-\int_{\Omega}|u_t|^{r-1}u_tudx+\|u\|_{p+1}^{p+1}.
  \end{split}
\end{equation*}
Adding and subtracting $(p+1)(1-\theta)E(t)$ with $\theta\in(0, 1)$ in the right hand side of the above equation, we have
\begin{equation}\label{3.2}
  \begin{split}\frac{d}{dt}\int_{\Omega}uu_tdx=&\frac{(p+1)(1-\theta)+2}{2}\|u_t\|_2^2+\frac{(p+1)(1-\theta)-2}{2}(\|\triangle u\|_2^2+\|\nabla u\|_2^2)\\[2mm]
  &+\frac{(p+1)(1-\theta)-2(\gamma+1)}{2(\gamma+1)}\beta\|\nabla u\|_2^{2(\gamma+1)}-\int_{\Omega}\nabla u\cdot \nabla u_tdx-\int_{\Omega}|u_t|^{r-1}u_tudx\\[2mm]
  &-(p+1)(1-\theta)E(t)+\theta\|u\|_{p+1}^{p+1}.
  \end{split}
\end{equation}
Using Young's inequality with $\epsilon\leq 1$, one gets 
\begin{equation}\label{3.3}
\big|\int_{\Omega}\nabla u\cdot \nabla u_tdx\big|\leq \frac{\epsilon}{2}\|\nabla u\|_2^2+\frac{1}{2\epsilon}\|\nabla u_t\|_2^2,
\end{equation}
\begin{equation}\label{3.4}
 \big|\int_{\Omega}|u_t|^{r-1}u_tudx\big|\leq \frac{r}{(r+1)\epsilon}\|u_t\|_{r+1}^{r+1}+\frac{\epsilon^r}{r+1}\|u\|_{r+1}^{r+1}.
\end{equation}
Using the interpolation inequality for $L^p-$norms, we obtain
\begin{equation}\label{3.5}
 \|u\|_{r+1}^{r+1}\leq s\|u\|_2^2+(1-s)\|u\|_{p+1}^{p+1}\,\, \text{with}\,\,\, s=\frac{p-r}{p-1}\in(0, 1].
\end{equation}
Inserting \eqref{3.3}-\eqref{3.5} into \eqref{3.2}, one can easily deduce that
\begin{equation}\label{3.6}
  \begin{split}\frac{d}{dt}\int_{\Omega}uu_tdx\geq&\frac{(p+1)(1-\theta)+2}{2}\|u_t\|_2^2+\frac{(p+1)(1-\theta)-2}{2}\|\triangle u\|_2^2-\frac{1}{2\epsilon}\|\nabla u_t\|_2^2\\[2mm]
  &+\frac{(p+1)(1-\theta)-2-\epsilon}{2}\|\nabla u\|_2^2+\frac{(p+1)(1-\theta)-2(\gamma+1)}{2(\gamma+1)}\beta\|\nabla u\|_2^{2(\gamma+1)}\\[2mm]
  &+\big(\theta -\frac{\epsilon^r(1-s)}{r+1}\big)\|u\|_{p+1}^{p+1}-\frac{\epsilon^rs}{r+1}\|u\|_2^2-\frac{r}{(r+1)\epsilon}\|u_t\|^{r+1}_{r+1}-(p+1)(1-\theta)E(t).
    \end{split}
\end{equation}
Now choosing $\theta=\frac{\epsilon^{r}(1-s)}{r+1}$ and $\epsilon \leq \delta_0:=\min\{1, [\frac{(p-2\gamma-1)(r+1)}{(p+1)(1-s)}]^{\frac{1}{r}}\}$ ($\delta_0=1$ when $r=1$, i.e. $s=1$), we have
$(p+1)(1-\theta)-2(\gamma+1)\geq0$ by $p>2\gamma+1$. Let
\begin{equation*}
  g(\epsilon):=(p+1)(1-\theta -2-\epsilon)=(p+1)(1-\frac{\epsilon^{r}(1-s)}{r+1})-2-\epsilon\,\,\text{for}\,\,\,\epsilon\in(0,\delta_0],
\end{equation*}
then we have
\begin{equation*}
  g'(\epsilon):=-\frac{(p+1)r\epsilon ^{r-1}(1-s)}{r+1}-1<0\,\,\text{for}\,\,\,\epsilon\in(0,\delta_0],
\end{equation*}
which implies that $g(\epsilon)$ is strictly decreasing in the interval $(0,\delta_0]$.

It follows from the fact $\displaystyle\lim_{\epsilon \rightarrow0}g(\epsilon)=p-1>0$ and the continuity of $g(\epsilon)$ that there exists $\delta_1>0$ such that
\begin{equation}\label{3.7}
g(\epsilon)>0\,\,\text{ for\,\, all}\,\, \epsilon\in (0, \delta_1).
\end{equation}
For $0<\epsilon <\delta_1<\delta_0$, using \eqref{2.5}, \eqref{3.7} and the embedding inequality $\|u\|_2\leq B_1\|\nabla u\|_2$, noticing the choice of $\theta$,  we can rewrite \eqref{3.6} as
\begin{equation}\label{3.8}
  \frac{d}{dt}\big(\int_{\Omega}uu_tdx-\frac{r}{(r+1)\epsilon}E(t)\big)\geq \frac{(p+1)(1-\theta)+2}{2}\|u_t\|_2^2+h(\epsilon)\|u\|_2^2-(p+1)(1-\theta)E(t),
\end{equation}
where
\begin{equation*}
  h(\epsilon)=\big[\frac{(p+1)(1-\frac{\epsilon^r(1-s)}{r+1})-2}{2}-\frac{\epsilon}{2}\big]\frac{1}{B_1^2}-\frac{\epsilon^r(1-s)}{r+1}\,\,\text{for}\,\,\epsilon\in(0, \delta_1)
\end{equation*}
By the similar argument as \eqref{3.7}, we can derive that  there exists $\delta_2\in(0, \delta_1)$ such that
\begin{equation}\label{3.9}
h(\epsilon)>0\,\,\text{ for\,\, all}\,\, \epsilon\in (0, \delta_2).
\end{equation}
By Cauchy-Schwarz inequality, \eqref{3.8} can be rewritten as
\begin{equation}\label{3.10}
\begin{split}
  \frac{d}{dt}\big(\int_{\Omega}uu_tdx-\frac{r}{(r+1)\epsilon}E(t)\big)\geq&   \sqrt{2[(p+1)(1-\theta)+2]h(\epsilon)}\int_{\Omega}uu_tdx-(p+1)(1-\theta)E(t)\\[2mm]
=&A(\epsilon)\bigg(\int_{\Omega}uu_tdx-B(\epsilon)E(t)\bigg),  \end{split}
\end{equation}
where
\begin{equation*}
  A(\epsilon)=\sqrt{2[(p+1)(1-\theta)+2]h(\epsilon)},
\end{equation*}
\begin{equation*}
  B(\epsilon)=\frac{(p+1)(1-\theta)}{A(\epsilon)}=\frac{(p+1)(1-\theta)}{\sqrt{2[(p+1)(1-\theta)+2]h(\epsilon)}}.
\end{equation*}
By the choice of $\theta$ and $h(\epsilon)$, we have
\begin{equation*}
  \lim_{\epsilon\rightarrow0^+}B(\epsilon)=\frac{(p+1)B_1}{\sqrt{(p+3)(p-1)}}<+\infty,\qquad   \lim_{\epsilon\rightarrow0^+}\frac{r}{(r+1)\epsilon}=+\infty,
\end{equation*}
which implies that there exists a sufficiently small $\delta_3\in (0, \delta_2)$ such that
\begin{equation*}
  B(\epsilon)\leq \frac{r}{(r+1)\epsilon}\quad \text{for\,any }\,\,\epsilon\in(0, \delta_3).
\end{equation*}
Hence, for any fixed sufficiently small $\epsilon_0\in(0, \delta_3)$, \eqref{3.10} can be rewritten as
\begin{equation}\label{3.11}
  \frac{d}{dt}\big(\int_{\Omega}uu_tdx-\frac{r}{(r+1)\epsilon_0}E(t)\big)\geq A(\epsilon_0)\big(\int_{\Omega}uu_tdx-\frac{r}{(r+1)\epsilon_0}E(t)\big),
\end{equation}
which implies  \eqref{3.1} holds with $A=A(\epsilon_0)$ and $B=\frac{r}{(r+1)\epsilon_0}$. From \eqref{3.11} and the above discussions, we can easily deduce that \eqref{3.1} also holds for the linear damping case ($r=1$).

\textbf{Step 2.} By contradiction, we suppose that $u$ is a global solution to the problem \eqref{1.1}. Recalling \eqref{3.1}, by Gronwall's inequality, we obtain
\begin{equation}\label{3.12}
  \int_{\Omega}uu_tdx-BE(t)\geq \big(\int_{\Omega}u_0u_1dx-BE(0)\big)e^{At}>0\,\,\text{for}\,\,t\geq0,
\end{equation}
where the assumption $0\leq E(0)<\frac{1}{B}\int_{\Omega}u_0u_1dx$ is used.  In view of $0\leq E(t)\leq E(0)$ and $\frac{d}{dt}\|u(t)\|_2^2=2\int_{\Omega}uu_tdx$,  from \eqref{3.12}, one has
\begin{equation}\label{3.13}
  \begin{split}\|u(t)\|_2^2=&\|u_0\|_2^2+2\int_0^t\int_{\Omega}uu_{\tau}dxd\tau\\[2mm]
  \geq &\|u_0\|_2^2+2\int_0^t\big(\int_{\Omega}u_0u_1dx-BE(0)\big)e^{A\tau}d\tau\\[2mm]
  =&\|u_0\|_2^2+ \frac{2}{A}(e^{At}-1)\big(\int_{\Omega}u_0u_1dx-BE(0)\big).
  \end{split}
\end{equation}

On the other hand, by \eqref{2.5} and H\"{o}lder's inequality, one has
\begin{equation*}
  \begin{split}\|u(t)\|_2=&\|u_0+\int_0^tu_t(\tau)d\tau\|_2\leq \|u_0\|_2+\int_{0}^t\|u_t(\tau)\|_2d\tau\\[2mm]
  \leq&\|u_0\|_2+|\Omega|^{\frac{r-1}{2(r+1)}}\int_{0}^t\|u_t(\tau)\|_{r+1}d\tau\\[2mm]
  \leq& \|u_0\|_2+|\Omega|^{\frac{r-1}{2(r+1)}}t^{\frac{r}{r+1}}\bigg(\int_{0}^t\|u_t(\tau)\|_{r+1}^{r+1}d\tau\bigg)^{\frac{1}{r+1}}\\[2mm]
  \leq& \|u_0\|_2+|\Omega|^{\frac{r-1}{2(r+1)}}t^{\frac{r}{r+1}}\bigg(E(0)-E(t)\bigg)^{\frac{1}{r+1}}\\[2mm]
  \leq& \|u_0\|_2+|\Omega|^{\frac{r-1}{2(r+1)}}t^{\frac{r}{r+1}}\big(E(0)\big)^{\frac{1}{r+1}},
  \end{split}
\end{equation*}
which is a contraction with \eqref{3.13} for $ t$ sufficiently large. Here we use the assumption  $u(t)$ is a global solution to the problem \eqref{1.1} and $E(t)\geq0$. Hence, $T_m<\infty$  and $u(t)$ blows up in finite time. This completes the proof.
\end{proof}

 From Theorem 3.1, we can show that the existence of finite time blow-up solutions with arbitrary initial energy level (including $E(0)>d$).
\begin{corollary}For any constant $R$ (including $R>d$), there exist two functions $u_0^R\in H$ and $u_1^R\in H_0^1(\Omega)$ satisfying $E(0)=R$, and the corresponding solution $u(t)$ blows up in finite time with initial data $u_0=u_0^R, \,u_1=u_1^R$, where $$E(0)=\frac{1}{2}\|u_1^R\|_2^2+\frac{1}{2}\|\triangle u_0^R\|_2^2+\frac{1}{2}\|\nabla u_0^R\|_2^2+\frac{1}{(2\gamma+1)}\beta\|\nabla u_0^R\|_2^{2(\gamma+1)}-\frac{1}{p+1}\|u_0^R\|_{p+1}^{p+1}.$$
\end{corollary}
\begin{proof}For any constant $R$,  choosing  two  arbitrary nonzero functions $v_1(x)$ and $v_2(x)$ such that $(v_1, v_2)=0$, we choose the initial data $u_0^R$ and $u_1^R$ as following
\begin{equation}\label{3.14}
 u_0^R(x):=r_1v_1(x),\,\,u_1^R(x):=r_1v_1(x)+r_2v_2(x),
\end{equation}
where $r_1$ and $r_2$ are two positive constants to be determined later. Let $u(t)$ denotes by the corresponding solution to the problem \eqref{1.1} with initial data $u_0=u_0^R, \,u_1=u_1^R$. Then the initial energy can be written as
\begin{equation}\label{3.15}
E(0)=\frac{r_2^2}{2}\|v_2\|_2^2+\chi(r_1),
\end{equation}
where
\begin{equation*}
  \chi(r_1)=\frac{r_1^2}{2}\|v_1\|_2^2+\frac{r_1^2}{2}\|\triangle v_1\|_2^2+\frac{r_1^2}{2}\|\nabla v_1\|_2^2+\frac{r_1^{2(\gamma+1)}}{(2\gamma+1)}\beta\|\nabla v_1\|_2^{2(\gamma+1)}-\frac{r_1^{p+1}}{p+1}\|v_1\|_{p+1}^{p+1}.
\end{equation*}
Noticing $1<2\gamma+1<p$, one has $\displaystyle\lim_{r_1\rightarrow \infty}\chi(r_1)=-\infty$. Hence, we can select $r_1$ sufficiently large such that $\chi(r_1)<R<\frac{1}{B}(u_0^R, u_1^R)=\frac{r_1^2}{B}\|v_1\|_2^2$. For such fixed $r_1$, choosing $r_2=\frac{\sqrt{2(R-\chi(r_1))}}{\|v_2\|_2}$, we have
\begin{equation*}
  E(0)=\frac{r_2^2}{2}\|v_2\|_2^2+\chi(r_1)=R.
\end{equation*}
Hence, we can derive from Theorem 3.1 that the corresponding solution $u(t)$ blows up in finite time with initial data $u_0=u_0^R, \,u_1=u_1^R$.

\end{proof}

Now, let us turn  our attention to the finite time blow up and an upper bound of the solution by adding some additional assumption to Theorem 3.1.
\begin{theorem} Let all the assumptions in Theorem 3.1 hold, and
\begin{equation}\label{3.16}
  \|u_0\|_2^2\geq \frac{p+2\gamma +3+\mu}{2\mu_0}E(0)\geq0,
\end{equation}where $\mu_0$ is defined by \eqref{3.19}, and $\mu$ is any positive number.  Then the solution $u(t)$ to the problem \eqref{1.1} blows up in finite time $T_m$ in the sense that $\displaystyle\lim_{t\rightarrow T_m^-}\|u\|_{p+1}=+\infty$ and an upper bound of $T_m$ can be estimated as follow
\begin{equation*}
  T_m\leq \frac{\mu_1}{\mu_2}\frac{1-\alpha}{\alpha}L^{-\frac{\alpha}{1-\alpha}}(0),
\end{equation*}
where $0<\displaystyle\alpha\leq\min\{\frac{p-1}{2(p+1)}, \,\frac{\gamma}{\gamma+1}\}$, $L(0)$, $\mu_1$ and $\mu_2$ will be determined in \eqref{3.27}, \eqref{3.21} and \eqref{3.26}, respectively.
\end{theorem}
\begin{proof} Without loss of generality, we can assume  $E(t)\geq0 $ for all $t\in [0, T_m)$.  Let $H(t)=E(0)-E(t)$ for any $t\geq0,$  then $H(t)\geq0.$  Define
\begin{equation*}
  L(t)=H^{1-\alpha}(t)+\varepsilon\bigg(\int_{\Omega}uu_tdx+\frac{1}{2}\|\nabla u\|_2^2\bigg),
\end{equation*}
where $\varepsilon$ will be determined later. We split the proof into the following three steps.

\textbf{Step 1.} Estimate  $L'(t)$. Due to \eqref{2.4} and the first equation of \eqref{1.1}, for any $\kappa\in (0, 1)$, by adding and subtracting $(p+1)(1-\kappa)\varepsilon E(t)$, one has
\begin{equation}\label{3.17}
\begin{split}L'(t)=&(1-\alpha)H^{-\alpha}(t)H'(t)+\varepsilon\|u_t\|_2^2-\varepsilon\|\triangle u\|_2^2-\varepsilon\|\nabla u\|_2^2-\varepsilon\beta\|\nabla u\|_2^{2(\gamma+1)}\\[2mm]
&-\varepsilon\int_{\Omega}|u_t|^{r-1}u_tudx+\varepsilon\|u\|_{p+1}^{p+1}\\[2mm]
=&(1-\alpha)H^{-\alpha}(t)H'(t)+\frac{(p+1)(1-\kappa)+2}{2}\varepsilon\|u_t\|_2^2+\frac{(p+1)(1-\kappa)-2}{2}\varepsilon\|\triangle u\|_2^2\\[2mm]
 &+\frac{(p+1)(1-\kappa)-2}{2}\varepsilon\|\nabla u\|_2^2+\frac{(p+1)(1-\kappa)-2(\gamma+1)}{2(\gamma+1)}\varepsilon\beta\|\nabla u\|_2^{2(\gamma+1)}\\[2mm]
 &+\varepsilon\kappa\|u\|_{p+1}^{p+1}-\varepsilon\int_{\Omega}|u_t|^{r-1}u_tudx-\varepsilon (p+1)(1-\kappa)E(t).
\end{split}
\end{equation}
To estimate $\int_{\Omega}|u_t|^{r-1}u_tudx$, by the similar argument as in Theorem 3.1, we  have
\begin{equation*}\begin{split}
 \big|\int_{\Omega}|u_t|^{r-1}u_tudx\big|\leq &\frac{r}{(r+1)\delta}\|u_t\|_{r+1}^{r+1}+\frac{\delta^r}{r+1}\|u\|_{r+1}^{r+1}\\[2mm]
 \leq&\frac{r}{(r+1)\delta}\|u_t\|_{r+1}^{r+1}+\frac{\delta^r}{r+1}\big(s\|u\|_2^2+(1-s)\|u\|_{p+1}^{p+1}\big).
\end{split}\end{equation*}
Let $\frac{r}{(r+1)\delta}=MH^{-\alpha}(t)$ for some positive constant $M$ to be chosen later. Noticing the assumption $E(t)\geq0$ implies $H(t)\leq E(0)$ and $$H'(t)=\|u_t\|_{r+1}^{r+1}+\|\nabla u_t\|_2^2,$$
one has
\begin{equation*}
  \frac{\delta^r}{r+1}=\frac{r^rH^{\alpha r}(t)}{(r+1)^{r+1}M^r}\leq\frac{r^rE^{\alpha r}(0)}{(r+1)^{r+1}M^r}.
\end{equation*}
Choosing $\kappa=\frac{p-(2\gamma+1)}{2(p+1)}\in(0,1)$ and inserting the above discussions into \eqref{3.17}, by \eqref{2.1}, we obtain
\begin{equation}\label{3.18}
\begin{split}L'(t)\geq&(1-\alpha-\varepsilon M)H^{-\alpha}(t)H'(t)+\frac{p+2\gamma+7}{4}\varepsilon\|u_t\|_2^2+\frac{p+2\gamma-1}{4}\varepsilon\|\nabla u\|_2^2\\[2mm]
 &+\frac{p-(2\gamma+1)}{2(\gamma+1)}\varepsilon\beta\|\nabla u\|_2^{2(\gamma+1)}+\varepsilon \big[\frac{p+2\gamma-1}{4}\lambda_1-\frac{r^rE^{\alpha r}(0)s}{(r+1)^{r+1}M^r}\big]\|u\|_2^2\\[2mm]
 &+\varepsilon\big[\frac{p-(2\gamma+1)}{2(p+1)}-\frac{r^rE^{\alpha r}(0)(1-s)}{(r+1)^{r+1}M^r}\big]\|u\|_{p+1}^{p+1}+\frac{p+2\gamma+3}{2}\varepsilon \big(H(t)-E(0)\big).
\end{split}
\end{equation}
Now, we can fix $M>0$ sufficiently large such that
\begin{equation}\label{3.19}
  \mu_0:=\frac{p+2\gamma-1}{4}\lambda_1-\frac{r^rE^{\alpha r}(0)s}{(r+1)^{r+1}M^r}>0,
\end{equation}
and
\begin{equation*}
  \zeta:=\frac{p-(2\gamma+1)}{2(p+1)}-\frac{r^rE^{\alpha r}(0)(1-s)}{(r+1)^{r+1}M^r}>0.
\end{equation*}
It follows from \eqref{3.13} and \eqref{3.16} that
\begin{equation*}
  \|u(t)\|_2^2\geq \|u_0\|_2^2\geq \frac{p+2\gamma +3+\mu}{2\mu_0}E(0).
\end{equation*}
Then, choosing $\varepsilon$ sufficiently small such that $1-\alpha-\varepsilon M\geq0$, we can rewrite \eqref{3.18} as
\begin{equation}\label{3.20}
\begin{split}L'(t)\geq&\varepsilon\big[\frac{p+2\gamma+7}{4}\|u_t\|_2^2+\frac{p+2\gamma-1}{4}\|\nabla u\|_2^2+\frac{p-(\gamma+1)}{2(\gamma+1)}\beta\|\nabla u\|_2^{2(\gamma+1)}\\[2mm]
&+\frac{\mu}{2}E(0)+\zeta\|u\|_{p+1}^{p+1}+\frac{p+2\gamma+3}{2} H(t)\big]\\[2mm]
\geq&\mu_1\big[\|u_t\|_2^2+\|\nabla u\|_2^2+\|\nabla u\|_2^{2(\gamma+1)}+\|u\|_{p+1}^{p+1}+H(t)+1\big],
\end{split}
\end{equation}
where $\mu_1$ is the smallest coefficient, that is
\begin{equation}\label{3.21}
  \mu_1=\varepsilon\min\{\frac{p+2\gamma-1}{4}, \frac{p-(2\gamma+1)}{2(\gamma+1)}\beta, \frac{p+2\gamma+3}{2}, \zeta, \frac{\mu}{2}E(0) \}.
\end{equation}

\textbf{Step 2.} Estimate  $L^{\frac{1}{1-\alpha}}(t)$.  In this step, we need to estimate
\begin{equation}\label{3.22}
  L^{\frac{1}{1-\alpha}}(t)=\big\{H^{1-\alpha}(t)+\varepsilon(\int_{\Omega}uu_tdx+\frac{1}{2}\|\nabla u\|_2^2)\big\}^{\frac{1}{1-\alpha}}.
\end{equation}
It is easy to derive the following inequality
\begin{equation*}
  \big|\int_{\Omega}uu_tdx\big|^{\frac{1}{1-\alpha}}\leq\|u_t\|_2^{\frac{1}{1-\alpha}}\|u\|_2^{\frac{1}{1-\alpha}}\leq\frac{1}{2(1-\alpha)}\|u_t\|_2^2+C_1\|u\|_{p+1}^{\frac{2}{2(1-\alpha)-1}}
\end{equation*}
where $C_1=\frac{2(1-\alpha)-1}{2(1-\alpha)}|\Omega|^{\frac{p-1}{p+1}\frac{1}{2(1-\alpha)-1}}>0$ ($0<\alpha<\frac{1}{2}$ by the choice of $\alpha$). Since $s_0:=\frac{2}{2(1-\alpha)-1}\leq p+1$, Young's inequality yields that
\begin{equation*}
 \|u\|_{p+1}^{s_0}\leq\frac{s_0}{(p+1)}\|u\|_{p+1}^{p+1}+\frac{p+1-s_0}{p+1}.
\end{equation*}
Hence, we obtain
\begin{equation}\label{3.23}
 \big|\int_{\Omega}uu_tdx\big|^{\frac{1}{1-\alpha}}\leq \frac{1}{2(1-\alpha)}\|u_t\|_2^2+\frac{C_1s_0}{(p+1)}\|u\|_{p+1}^{p+1}+\frac{p+1-s_0}{p+1}C_1.
\end{equation}
Since $0<\alpha \leq\frac{\gamma}{\gamma+1}$, one has the following estimate (see \cite{LiaoSCM2023,LiaoLiMJM2023} for details)
\begin{equation}\label{3.24}
 \|\nabla u\|_2^{\frac{2}{1-\alpha}}\leq C_2\big(\|\nabla u\|_2^2+\|\nabla u\|_2^{2(\gamma+1)}\big)
\end{equation}
for some positive constant $C_2$.
Combining \eqref{3.22}-\eqref{3.24} with the following   algebraic inequality
\begin{equation*}
  (a+b+c)^{l}\leq 2^{2(l-1)}(a^l+b^l+c^l), \,\text{for}\,a, b, c\geq0\,\text{and}\,l\geq1,
\end{equation*}
 one has
\begin{equation}\label{3.25}
\begin{split}L^{\frac{1}{1-\alpha}}(t)\leq&2^{\frac{2\alpha}{1-\alpha}}\bigg(H(t)+\varepsilon^{\frac{1}{1-\alpha}}|\int_{\Omega}uu_tdx|^{\frac{1}{1-\alpha}}+(\frac{\varepsilon}{2})^{\frac{1}{1-\alpha}}\|\nabla u\|_2^{\frac{2}{1-\alpha}}\bigg)\\[2mm]
\leq&\mu_2\big[\|u_t\|_2^2+\|\nabla u\|_2^2+\|\nabla u\|_2^{2(\gamma+1)}+\|u\|_{p+1}^{p+1}+H(t)+1\big],
\end{split}
\end{equation}
where
\begin{equation}\label{3.26}
  \mu_2=2^{\frac{2\alpha}{1-\alpha}}\max\{\frac{\varepsilon^{\frac{1}{1-\alpha}}}{2(1-\alpha)}, \varepsilon^{\frac{1}{1-\alpha}}\frac{C_1s_0}{(p+1)},   \varepsilon^{\frac{1}{1-\alpha}}\frac{p+1-s_0}{p+1}C_1, (\frac{\varepsilon}{2})^{\frac{1}{1-\alpha}}C_2 \}.
\end{equation}

\textbf{Step 3.} Complete the proof.  It follows from \eqref{3.20} and \eqref{3.25} that $L'(t)\geq \frac{\mu_1}{\mu_2}L^{\frac{1}{1-\alpha}}(t)$. By the assumption, we obtain
 \begin{equation}\label{3.27}
  L(0)=\varepsilon(\int_{\Omega}u_0u_1dx+\frac{1}{2}\|\nabla u_0\|_2^2)>0
\end{equation}
 Hence we have
\begin{equation*}
  L^{\frac{\alpha}{1-\alpha}}(t)\geq \frac{1}{L^{-\frac{\alpha}{1-\alpha}}(0)-\frac{\mu_2\alpha}{\mu_1(1-\alpha)}t},
\end{equation*}
which implies that there exists a finite time $T_m$ such that $$\displaystyle\lim_{t\rightarrow T_m^-}L(t)=+\infty \, \,\text{and}\,\,  T_m\leq \frac{\mu_1}{\mu_2}\frac{1-\alpha}{\alpha}L^{-\frac{\alpha}{1-\alpha}}(0)$$.

Finally, we show that $\displaystyle\lim_{t\rightarrow T_m^-}L(t)=+\infty$ implies $\displaystyle\lim_{t\rightarrow T_m^-}\|u\|_{p+1}=+\infty$. We divide the proof into three cases due to the definition of $L(t)$.\\
(a) $H(t)\rightarrow+\infty$. This case is impossible by the assumption $E(t)\geq0$ for all $t\in[0,T_m)$.\\
(b) If $\displaystyle\int_{\Omega}uu_tdx\rightarrow+\infty$, it follows from \eqref{2.1} and Cauchy's inequality that 
\begin{equation*}
  \int_{\Omega}uu_tdx\leq \frac{1}{2}\|u_t\|_2^2+\frac{1}{2}\|u\|_2^2\leq \frac{1}{2}\|u_t\|_2^2+\frac{1}{2\lambda_1}\|\triangle u\|_2^2.
\end{equation*}
Noticing $E(t)\leq E(0)$ and \eqref{2.4}, we have
\begin{equation}\label{3.28}\begin{split}
  \frac{1}{2}(\|u_t\|_2^2+\|\triangle u\|_2^2)=&E(t)+\frac{1}{p+1}\|u\|_{p+1}^{p+1}-\frac{1}{2}\|\nabla u\|_2^2-\frac{1}{2(\gamma+1)}\beta\|\nabla u\|_{2}^{2(\gamma +1)}\\[2mm]
  \leq& E(0)+\frac{1}{p+1}\|u\|_{p+1}^{p+1}.
\end{split}\end{equation}
Combining  the above two  inequalities, we easily have $\displaystyle\lim_{t\rightarrow T_m^-}\|u\|_{p+1}=+\infty$  from the fact  $\displaystyle\int_{\Omega}uu_tdx\rightarrow+\infty$ as $t\rightarrow T_m^-$.\\
(c) If $\|\nabla u\|_2^2\rightarrow +\infty$, in view of \eqref{3.28} and \eqref{2.1},  we can also derive  $\displaystyle\lim_{t\rightarrow T_m^-}\|u\|_{p+1}=+\infty$.

\end{proof}

Now we shall show the assumption $E(t)\geq0$ for all $t\in [0, T_m)$ is valid. In fact, we can also derive the finite time blow-up result and an upper bound of the blow-up time when there exist $t_0\in[0, T_m)$ such that $E(t_0)<0$.
 Without loss of generality, we can assume that $E(0)<0 $ for convenience.
\begin{theorem} Let $u_0\in H$, $u_1\in H_0^1(\Omega)$ and \eqref{1.2} hold. Assume that the initial energy $E(0)<0$, then the solution $u(t)$ to the problem \eqref{1.1} blows up in finite time $T_m$ in the sense that $\displaystyle\lim_{t\rightarrow T_m^-}\|u\|_{p+1}=+\infty$ and an upper bound of $T_m$ can be estimated as follow
\begin{equation*}
  T_m\leq \frac{\mu_3}{\mu_4}\frac{1-\alpha}{\alpha}L^{-\frac{\alpha}{1-\alpha}}(0),
\end{equation*}
where $\alpha$, $L(0)>0$, $\mu_3$ and $\mu_4$ will be determined in \eqref{3.30}, \eqref{3.32}, \eqref{3.36} and \eqref{3.38}, respectively.
\end{theorem}
\begin{proof} The proof is similar as that of Theorem 3.1. Since there is no assumption \eqref{3.16} and lack the monotone increasing property of norm $\|u\|_2^2$ (see\eqref{3.13}), we must modify some estimates in Theorem 3.2.  Using the same notations in Theorem 3.2, we give a sketch of the proof.

Let $H(t)=-E(t)$ and
\begin{equation}\label{3.29}
  L(t)=H^{1-\alpha}(t)+\varepsilon\bigg(\int_{\Omega}uu_tdx+\frac{1}{2}\|\nabla u\|_2^2\bigg),
\end{equation}where
\begin{equation}\label{3.30}
  0<\displaystyle\alpha\leq\min\{\frac{p-r}{(p+1)r},\frac{p-1}{2(p+1)}, \,\frac{\gamma}{\gamma+1}\}.
\end{equation}
Noticing  \eqref{2.4} and the definition of $H(t)$,  we have
\begin{equation}\label{3.31}
  0<H(0)\leq H(t)\leq \frac{1}{p+1}\|u\|_{p+1}^{p+1}.
\end{equation}
By \eqref{3.29} and \eqref{3.31}, we can choose sufficiently small $\varepsilon$ such that
\begin{equation}\label{3.32}
  L(0)=(H(0))^{1-\alpha}+\varepsilon\big(\int_{\Omega}u_0u_1dx+\frac{1}{2}\|\nabla u_0\|_2^2\big)>0.
\end{equation}
  Since $p>r$, using the following  algebraic inequality (see \cite{Chen2012})
$$
z^\nu \leq(z+1) \leq\left(1+\frac{1}{a}\right)(z+a) \quad \text { for all } z \geq 0,0 \leq \nu \leq 1, a > 0,
$$
 for any $s\in [2, p+1]$, we obtain
\begin{equation}\label{3.33}
 \|u\|_{p+1}^{s} \leq\left(1+\frac{1}{H(0)}\right)\left(\|u\|_{p+1}^{p+1}+H(0)\right) \leq \left(1+\frac{1}{H(0)}\right)\left(\|u\|_{p+1}^{p+1}+H(t)\right).
\end{equation}
Hence, following H\"{o}lder's inequality, Young's inequality and \eqref{3.31}, one has
\begin{equation*}\begin{split}
 \big|\int_{\Omega}|u_t|^{r-1}u_tudx\big|\leq &\frac{rH^{-\alpha}(t)}{(r+1)\delta}\|u_t\|_{r+1}^{r+1}+\frac{\delta^rH^{\alpha r}(t)}{r+1}\|u\|_{r+1}^{r+1}\\[2mm]
 \leq&\frac{rH^{-\alpha}(t)}{(r+1)\delta}\|u_t\|_{r+1}^{r+1}+\frac{\delta^r}{r+1}(\frac{1}{p+1})^{\alpha r}\|u\|_{r+1}^{(p+1)\alpha r+r+1}\\[2mm]
 \leq&\frac{rH^{-\alpha}(t)}{(r+1)\delta}\|u_t\|_{r+1}^{r+1}+C_3\delta^r\big(H(t)+\|u\|_{p+1}^{p+1}\big),
\end{split}\end{equation*}
where $C_3=(1+\frac{1}{H(0)})|\Omega|^{\frac{p-r-(p+1)\alpha r}{p+1}}\frac{1}{r+1}(\frac{1}{p+1})^{\alpha r}.$  Here we also use \eqref{3.30}, that is $(p+1)\alpha r+r+1\leq p+1$. Adding and subtracting $\varepsilon(p+1)(1-\kappa)E(t)$ with the same choice $\kappa=\frac{p-(2\gamma+1)}{2(p+1)}$ as that in the proof of Theorem 3.2,  we can rewrite \eqref{3.18} as 
\begin{equation}\label{3.34}
\begin{split}L'(t)\geq&\big(1-\alpha-\varepsilon \frac{r}{(r+1)\delta}\big)H^{-\alpha}(t)H'(t)+\frac{p+2\gamma+7}{4}\varepsilon\|u_t\|_2^2\\[2mm]
 &+\frac{p+2\gamma-1}{4}\varepsilon(\|\nabla u\|_2^2+\|\triangle u\|_2^2)+\frac{p-(2\gamma+1)}{2(\gamma+1)}\varepsilon\beta\|\nabla u\|_2^{2(\gamma+1)} \\[2mm]
 &+\varepsilon\big[\frac{p-(2\gamma+1)}{2(p+1)}-C_3\delta^r\big]\|u\|_{p+1}^{p+1}+\varepsilon\big[\frac{p+2\gamma+3}{2}-C_3\delta^r\big]H(t) .
\end{split}
\end{equation}
Now, we can fix $\delta>0$ sufficiently small such that
\begin{equation*}
  \frac{p-(2\gamma+1)}{2(p+1)}-C_3\delta^r>0\,\, \text{and}\,\, \frac{p+2\gamma+3}{2}-C_3\delta^r>0.
\end{equation*}
We can select $\varepsilon$ sufficient small such that $1-\alpha-\varepsilon \frac{r}{(r+1)\delta}\geq0$. Then \eqref{3.34} yields that
\begin{equation}\label{3.35}
L'(t)\geq \mu_3\big[\|u_t\|_2^2+\|\nabla u\|_2^2+\|\triangle u\|_2^2+\|\nabla u\|_2^{2(\gamma+1)}+\|u\|_{p+1}^{p+1}+H(t)\big].
\end{equation}
where \begin{equation}\label{3.36}
  \mu_3=\varepsilon\min\{\frac{p+2\gamma-1}{4}, \frac{p-(2\gamma+1)}{2(\gamma+1)}\beta,  \frac{p-(2\gamma+1)}{2(p+1)}-C_3\delta^r, \frac{p+2\gamma+3}{2}-C_3\delta^r \}.
\end{equation}

On the other hand, by \eqref{3.33}, we can modify  \eqref{3.23} as follow
\begin{equation*}\begin{split}
  \big|\int_{\Omega}uu_tdx\big|^{\frac{1}{1-\alpha}}\leq&\frac{1}{2(1-\alpha)}\|u_t\|_2^2+C_1\|u\|_{p+1}^{\frac{2}{2(1-\alpha)-1}}\\[2mm]
  \leq& \frac{1}{2(1-\alpha)}\|u_t\|_2^2+C_1\big(1+\frac{1}{H(0)}\big)\big(H(t)+\|u\|_{p+1}^{p+1}\big).
\end{split}\end{equation*}
Then, by the similar argument as  Step 2 of Theorem 3.2, we can rewrite \eqref{3.25} as
\begin{equation}\label{3.37}
  L^{\frac{1}{1-\alpha}}(t)\leq \mu_4\big[\|u_t\|_2^2+\|\nabla u\|_2^2+\|\nabla u\|_2^{2(\gamma+1)}+\|u\|_{p+1}^{p+1}+H(t)\big],
\end{equation}

where
\begin{equation}\label{3.38}
  \mu_4=2^{\frac{2\alpha}{1-\alpha}}\max\{\frac{\varepsilon^{\frac{1}{1-\alpha}}}{2(1-\alpha)}, 1+\varepsilon^{\frac{1}{1-\alpha}}\frac{C_1(1+\frac{1}{H(0)})}{(p+1)},   (\frac{\varepsilon}{2})^{\frac{1}{1-\alpha}}C_2 \}.
\end{equation}
 It follows from \eqref{3.35} and \eqref{3.37} that $L'(t)\geq \frac{\mu_3}{\mu_4}L^{\frac{1}{1-\alpha}}(t)$, which implies that $\displaystyle\lim_{t\rightarrow T_m^-}L(t)=+\infty$. Meantime, the blow-up time $T_m$ can be estimated from
above as $  T_m\leq \frac{\mu_3}{\mu_4}\frac{1-\alpha}{\alpha}L^{-\frac{\alpha}{1-\alpha}}(0)$.

 Now we remain to prove $\displaystyle\lim_{t\rightarrow T_m^-}L(t)=+\infty$ implies $\displaystyle\lim_{t\rightarrow T_m^-}\|u\|_{p+1}=+\infty$. Indeed, the  proof is the same as that in Theorem 3.2 except Case (a).
If $H(t)\rightarrow+\infty$, one can easily obtain $\displaystyle\lim_{t\rightarrow T_m^-}\|u\|_{p+1}=+\infty$ from \eqref{3.31}.

\end{proof}

Finally, we will  give  lower bounds of the blow-up time.
\begin{theorem} Let $u_0\in H$, $u_1\in H_0^1(\Omega)$ and $1<p\leq\frac{n}{n-4}$ ($n\geq 5$) hold.  Let $u(t)$ be the solution to the problem \eqref{1.1}, which blows up at  a finite time $T$, then the lower bound for the blow-up time $T$ can be given by
\begin{equation*}
  T_m\geq \int_{F(0)}^{+\infty}\frac{1}{K_1+y+K_2y^p}dy
\end{equation*}
where $F(0)=\|u_0\|_{p+1}^{p+1}$, $K_1$ and $K_2$ are defined by \eqref{3.41}.
\end{theorem}
\begin{proof}It follows from \eqref{2.4} and \eqref{2.5} that $E(t)\leq E(0):=\varpi$ and
\begin{equation}\label{3.39}\begin{split}
  \frac{1}{2}\|u_t\|_2^2+\frac{1}{2}\|\nabla u\|_2^2+\frac{1}{2}\|\triangle u\|_2^2+\frac{\beta}{2(\gamma+1)}\|\nabla u\|_2^{2\gamma +2}&=E(t)+ \frac{1}{p+1}\|u\|_{p+1}^{p+1} \\[2mm]
  &\leq\varpi+\frac{1}{p+1}\|u\|_{p+1}^{p+1}.
\end{split}\end{equation}
Now let us denote by $F(t)=\|u(t)\|_{p+1}^{p+1}$. Noticing $1<p\leq\frac{n}{n-4}$, let $B_*$ be the best embedding constant $H_0^2(\Omega)\hookrightarrow L^{2p}$, then by Young's inequality and \eqref{3.39}, one has
\begin{equation}\label{3.40}
  \begin{split}
  F'(t)=&(p+1)\int_{\Omega}u^{p-1}uu_tdx \leq (p+1)\big(\frac{1}{2}\|u_t\|_2^2+\frac{1}{2}\|u\|_{2p}^{2p}\big)\\[2mm]
  \leq& (p+1)\big(\varpi+\frac{1}{p+1}F(t)+\frac{B_*^{2p}}{2}\|\triangle u\|_2^{2p}\big)\\[2mm]
  \leq& (p+1)\bigg(\varpi+\frac{1}{p+1}F(t)+\frac{B_*^{2p}}{2}(2\varpi+\frac{2}{p+1}F(t))^p\bigg)\\[2mm]
  \leq& (p+1)\bigg(\varpi+\frac{1}{p+1}F(t)+B_*^{2p}2^{p-2}\big((2\varpi)^{p}+(\frac{2}{p+1}F(t))^p\big)\bigg)\\[2mm]
  =&K_1+F(t)+K_2F^p(t),
  \end{split}\end{equation}
where
\begin{equation}\label{3.41}
K_1=(p+1)(\varpi+B_*^{2p}2^{p-2}\big((2\varpi)^{p})\,\, \text{and}\, K_2=B_*^{2p}2^{2p-2}(p+1)^{-p}.
\end{equation}
Noticing $\displaystyle\lim_{t\rightarrow T_m}\|u\|_{p+1}=+\infty$, we obtain  from \eqref{3.40} that
\begin{equation*}
  T_m\geq \int_{F(0)}^{+\infty}\frac{1}{K_1+y+K_2y^p}dy.
\end{equation*}
Hence, we complete the proof.
\end{proof}

 By making full use of the strong damping term $-\triangle u_t$, we will  give another lower bound of the blow-up time under the condition \eqref{1.2}.
\begin{theorem} Let $u_0\in H$, $u_1\in H_0^1(\Omega)$ and \eqref{1.2} hold. Assume that $u(t)$ is the solution of the problem \eqref{1.1}, which blows up at  a finite time $T_m$, then the lower bound for the blow-up time $T_m$  can be given as
$T_m\geq K_3G^{1-p}(0)$, where $G(0)$ and $K_3$ are defined in \eqref{3.43}.
\end{theorem}
\begin{proof} Define
\begin{equation*}
  G(t):=\frac{1}{2}\|u_t\|_2^2+\frac{1}{2}\|\nabla u\|_2^2+\frac{1}{2}\|\triangle u\|_2^2+\frac{\beta}{2(\gamma+1)}\|\nabla u\|_2^{2\gamma +2}=E(t)+\frac{1}{p+1}\|u\|_{p+1}^{p+1}
\end{equation*}
In view of \eqref{2.5}, we have
\begin{equation*}
  G'(t)=E'(t)+\int_{\Omega}|u|^{p-1}uu_tdx\leq -\|\nabla u_t\|_2^2+\int_{\Omega}|u|^{p-1}uu_tdx.
\end{equation*}
By H\"{o}lder's inequality, Young's inequality, we obtain
\begin{equation}\label{3.42}\begin{split}
 G'(t)&\leq -\|\nabla u_t\|_2^2+ \|u_t\|_{\frac{2n}{n-2}}\|u\|^p_{\frac{2np}{n+2}}\leq -\|\nabla u_t\|_2^2+C\|\nabla u_t\|_2\|\triangle u\|_2^p\\[2mm]
 &\leq C\|\triangle u\|_2^{2p}\leq CG^p(t),
\end{split}\end{equation}
for some positive constant $C$,  where we use the embedding $\|u_t\|_{\frac{2n}{n-2}}\leq C_*\|\nabla u_t\|$ and \\$ \|u\|_{\frac{2np}{n+2}}\leq C^*\|\triangle u\|_2,$ since $\frac{2np}{n+2}\leq \frac{2n}{n-4}$.

Noticing $\displaystyle\lim_{t\rightarrow T_m^-}\|u\|_{p+1}=+\infty$ implies $\displaystyle\lim_{t\rightarrow T_m^-}G(t)=+\infty$, integrating inequality \eqref{3.42} with $t$ over $(0, T_m)$ we have
\begin{equation*}
  T_m\geq \int_{0}^{T_m}C^{-1}G^{-p}(t)G'(t)dt=\int_{G(0)}^{+\infty}C^{-1}y^{-p}dy=K_3G^{1-p}(0),
\end{equation*}
where
\begin{equation}\label{3.43}
   G(0)=\frac{1}{2}\|u_1\|_2^2+\frac{1}{2}\|\nabla u_0\|_2^2+\frac{1}{2}\|\triangle u_0\|_2^2+\frac{\beta}{2(\gamma+1)}\|\nabla u_0\|_2^{2\gamma +2}>0\,\,\text{and} \,\,K_3=\frac{G^{1-p}(0)}{(p-1)C}>0.
\end{equation}

\end{proof}

\section*{Acknowledgments}
The authors would like to thank the referees for the careful reading
of this paper and for the valuable suggestions to improve the
presentation and the style of the paper. This paper is supported by  the
Innovative Funds Plan of Henan University of Technology (No.2020ZKCJ09) and National Natural Science Foundation of China (No.11801145).


\end{document}